\documentstyle[11pt,twoside]{article}

\newtheorem{theorem}{\noindent\bf Theorem}[section]
\newtheorem{proposition}[theorem]{\noindent\bf Proposition}
\newtheorem{lemma}[theorem]{\noindent\bf Lemma}
\newtheorem{corollary}[theorem]{\noindent\bf Corollary}

\newtheorem{remark}[theorem]{\noindent\bf Remark}

\newtheorem{example}[theorem]{\noindent\bf Example}

\newtheorem{definition}[theorem]{\noindent\bf Definition}

\newcommand{\e}{\hfill\blacksquare}
\input{amssym}
\begin{document}
\date{}
\title{{\Large\bf Amalgamated duplication of the
Banach algebra $\bf{\frak A}$\\ along a ${\frak A}$-bimodule
${\mathcal A}$}}
\author{{\normalsize\sc Hossein Javanshiri and Mehdi Nemati}}

\maketitle
\normalsize

\begin{abstract}
Let ${\mathcal A}$ and ${\frak A}$ be Banach algebras such that
${\mathcal A}$ is a Banach ${\frak A}$-bimodule with compatible
actions. We define the product ${\cal A}\rtimes{\frak A}$, which
is a strongly splitting Banach algebra extension of ${\frak A}$ by
$\cal A$. After characterization of the multiplier algebra,
topological centre, (maximal) ideals and spectrum of ${\cal
A}\rtimes{\frak A}$, we restrict our investigation to the study of
semisimplicity, regularity, Arens regularity of ${\cal
A}\rtimes{\frak A}$ in relation to that of the algebras $\cal A$,
$\frak A$ and the action of $\frak A$ on $\cal A$. We also compute
the first cohomology group $H^1{(}{\cal A}\rtimes{\frak A},({\cal
A}\rtimes{\frak A})^{(n)}{)}$ for all $n\in {\Bbb N}\cup\{0\}$ as
well as the first-order cyclic cohomology group
$H_\lambda^1{(}{\cal A}\rtimes{\frak A},({\cal A}\rtimes{\frak
A})^{(1)}{)}$, where $({\cal A}\rtimes{\frak A})^{(n)}$ is the
n-th dual space of ${\cal A}\rtimes{\frak A}$ when $n\in{\Bbb N}$
and ${\cal A}\rtimes{\frak A}$ itself when $n=0$. These results
are not only of interest in their own right, but also they pave
the way for obtaining some new results for Lau products and module
extensions of Banach algebras as well as triangular Banach
algebra. Finally, special attention is devoted to the cyclic and
$n$-weak amenability of ${\cal A}\rtimes{\frak A}$.\\


\end{abstract}


\section{Introduction}

Let $\cal A$ and $\cal B$ be Banach algebras. It is well-known
that the Cartesian product space ${\cal A}\times {\cal B}$
equipped with the $\ell^1$-norm and coordinatewise operations is a
Banach algebra. In order to provide new examples of Banach
algebras as well as a wealth of (counter) examples in different
branches of functional analysis, the construction of an algebra
product on the Cartesian product space ${\cal A}\times {\cal B}$
to make it a Banach algebra has been studied in a series of papers
recently. The first important paper related to this construction
is Lau's paper \cite{Lau} which he defined a new algebra product
on the Cartesian product space ${\cal A}\times {\cal B}$ for the
case where $\cal B$ is the predual of a van Neumann algebra $\cal
M$ such that the identity of $\cal M$ is a multiplicative linear
functional on $\cal B$. Later on, Monfared \cite{M} extended the
Lau product to arbitrary Banach algebras. This construction has
many applications in different contexts, see for examples \cite{
vishki,kaniuth,Mon2,nemati1}. Moreover, it is notable that such
constructions are also investigated in commutative ring theory and
have been extensively studied in recent years, see for examples
\cite{ring0,ring1,ring2,ring3} and the references therein.

On the other hand, motivated by Wedderburn's principal theorem
\cite{bad2}, splitting of Banach algebra extensions has been a
major question in the theory of Banach algebras and several
researchers have used the splitting of Banach algebra extensions
as a tool for the study of Banach algebras. For example, module
extensions as generalizations of Banach algebra extensions were
introduced and first studied by Gourdeau \cite{gro} were used to
show that amenability of ${\cal A}^{(2)}$, the second dual space
of $\cal A$, implies amenability of $\cal A$; Filali and Eshaghi
Gordji \cite{esfi} used triangular Banach algebras to answer some
questions asked by Lau and \"{U}lger \cite{lauul}; Zhang
\cite{zhang} used module extensions to construct an example of a
weakly amenable Banach algebra which is not 3-weakly amenable. For
some other applications of splitting of Banach algebra extensions
we refer the reader to the references \cite{bd,bdl,ob,fel,thomas}.

This paper continues these investigations. In details, the product
${\cal A}\rtimes{\frak A}$ of Banach algebras defined in this
paper (Definition \ref{defmain} below) and we point out that the
Lau products and module extension of Banach algebras as well as
triangular Banach algebras can be regarded as examples of ${\cal
A}\rtimes{\frak A}$ (Example \ref{exmain} below) and many basic
properties can be derived within this more general context. Apart
from the determination of the action of ${\cal A}\rtimes{\frak A}$
on its n-th dual space $({\cal A}\rtimes{\frak A})^{(n)}$ ($n\in
{\Bbb N}$), we give a characterization of the spectrum of ${\cal
A}\rtimes{\frak A}$. Sufficient and necessary conditions for
semisimplicity, regularity, Arens regularity and strong Arens
irregularity of ${\cal A}\rtimes{\frak A}$ are given. Moreover, we
obtain characterizations of multiplier algebra, topological
centre, (maximal) ideals of ${\cal A}\rtimes{\frak A}$ and compute
the first cohomology group $H^1{(}{\cal A}\rtimes{\frak A},({\cal
A}\rtimes{\frak A})^{(n)}{)}$ ($n\in {\Bbb N}\cup\{0\}$) as well
as the first-order cyclic cohomology group $H_\lambda^1{(}{\cal
A}\rtimes{\frak A},({\cal A}\rtimes{\frak A})^{(1)}{)}$.


\section{\large\bf Preliminaries}

Throughout this paper, for any Banach space $\cal Y$ and integer number $n\geq 0$, by ${\cal Y}^{(n)}$ we denote the $n$-th dual space of $\cal Y$ when $n\geq 1$ and $\cal Y$ itself when $n=0$. Moreover, we always consider $\cal Y$ as naturally embedded into ${\cal Y}^{(2)}$.

Let $\cal C$ be a Banach algebra and let $\sigma({\cal C})$ stands for the spectrum of $\cal C$; that is, the set of all non-zero multiplicative linear functionals on $\cal C$. If ${\mathcal C}$ is commutative
and $c\in{\cal C}$, then we define $\hat{c}:\sigma({\cal C})\rightarrow{\Bbb C}$ by $\hat{c}(\varphi)=\varphi(c)$ for all $\varphi\in \sigma({\mathcal C})$. The continuous function
$\hat{c}$ is called the Gelfand transform of $c$. It is
easily checked that the mapping
$$\Gamma_{\cal C}:{\cal C}\rightarrow C_0(\sigma({\cal C}));\quad c\mapsto\hat{c}$$
is a homomorphism which is called the Gelfand homomorphism of
${\cal C}$. In this  case, we say that $\cal C$ is semisimple if the Gelfand homomorphism is injective and it is called regular if for each
closed subset $E$ of $\sigma({\cal C})$ and $\varphi_0\in\sigma({\cal C})\setminus E$, there exists $c\in{\cal C}$ such that $\varphi_0(c)\neq 0$ and $\varphi(c)=0$ for
all $\varphi\in E$. Recall also that on ${\cal C}^{(2)}$ there exist two natural products
extending the one on $\cal C$, known as the first and second Arens products of ${\cal C}^{(2)}$ where the first Arens product ``$\circ$" on ${\cal C}^{(2)}$ defined in three steps as follows:
\begin{eqnarray*}
\big<a^{(2)}\circ b^{(2)},a^{(1)}\big>&=&\big<a^{(2)},b^{(2)}\circ a^{(1)}\big>,\\
\big<b^{(2)}\circ a^{(1)},a\big>&=&\big<b^{(2)},a^{(1)}\circ a\big>,\\
\big<a^{(1)}\circ a,b\big>&=&\big<a^{(1)},ab\big>,
\end{eqnarray*}
where $a,b\in{\cal C}$, $a^{(1)}\in{\cal C}^{(1)}$ and $a^{(2)}, b^{(2)}\in {\cal C}^{(2)}$; similarly, by using symmetry, the second Arens product ``$\vartriangle$" on ${\cal C}^{(2)}$ is defined as follows:
\begin{eqnarray*}
\big<a^{(2)}\vartriangle b^{(2)},a^{(1)}\big>&=&\big<b^{(2)},a^{(1)}\vartriangle a^{(2)}\big>,\\
\big<a^{(1)}\vartriangle a^{(2)},a\big>&=&\big<a^{(2)},a\vartriangle a^{(1)}\big>,\\
\big<a\vartriangle a^{(1)},b\big>&=&\big<a^{(1)},ba\big>,
\end{eqnarray*}
The first topological centre
$${\frak Z}_t({\cal C}^{(2)}):={\Big\{}a^{(2)}\in{\cal C}^{(2)}:~a^{(2)}\circ b^{(2)}=a^{(2)}\vartriangle b^{(2)}~{\hbox{for all}}~b^{(2)}\in{\cal C}^{(2)}{\Big\}}$$
of ${\cal C}^{(2)}$ is a closed subalgebra of $({\cal
C}^{(2)},\circ)$ containing $\cal C$. The algebra $\cal C$ is
called Arens regular [respectively, strongly Arens irregular] if
${\frak Z}_t({\cal C}^{(2)})={\cal C}^{(2)}$   [respectively,
${\frak Z}_t({\cal C}^{(2)})= {\cal C}$].

Now, let ${\cal X}$ be a Banach $\cal C$-bimodule. Then
${\cal X}^{(1)}$ is also
a Banach $\cal C$-bimodule by the following module actions:
$$\langle x^{(1)}\cdot c,x \rangle=\langle x^{(1)},c\cdot x
\rangle,\quad\quad\quad \langle c\cdot x^{(1)},x \rangle=\langle
x^{(1)},x\cdot c \rangle,$$ where $c\in{\cal C},~x\in {\cal
{\cal X}},~x^{(1)}\in {\cal {\cal X}}^{(1)}$. We may apply the same argument to show that
for each $n\in {\Bbb N}$ ${\cal X}^{(n)}$  is a Banach $\cal C$-bimodule.
A {derivation} from $\cal C$ into ${\cal X}$ is a
continuous linear map $D:{\cal C}\rightarrow {\cal X}$
such that for every $c_1,c_2 \in{\cal C}$
$$D(c_1c_2)=D(c_1)\cdot c_2+c_1\cdot D(c_2).$$ A derivation
$D:{\cal C}\rightarrow{\cal C}^{(1)}$ is called cyclic if
$$\big<D(c_1),c_2\big>+\big<D(c_2),c_1\big>=0\quad\quad\quad\quad\quad\quad(c_1,c_2\in{\cal C}),$$
For $x\in {\cal X}$, define ${\verb"ad"}_{x}$ from $\cal C$ into
${\cal X}$ by ${\verb"ad"}_{x}(c)=c\cdot x-x\cdot c$. It is easy
to show that ${\verb"ad"}_{x}$ is a derivation; such derivations
are called {inner derivations}. We denote the set of all
derivations from ${\cal C}$ into ${\cal X}$ by $Z^1({\cal C},{\cal
X})$. The {first cohomology group} $H^1({\cal C}, {\cal X})$ is
the quotient of the space of all  derivations by the space of all
inner derivations, and similarly the {first-order cyclic
cohomology group} $H_\lambda^1{(}{\cal C},{\cal X}{)}$ is the
quotient of the space of all  cyclic derivations by the space of
all inner derivations. In many situations triviality of these
spaces are of considerable importance. In particular, $\cal C$ is
called {contractible} [respectively, {amenable}] if $H^1({\cal C}
,{\cal X})=0$ [respectively, $H^1({\cal C} ,{\cal X}^{(1)})=0$]
for every Banach ${\cal C}$-bimodule ${\cal X}$ and for $n\in
{\Bbb N}\cup\{0\}$ it is called {$n$-weakly amenable} if
$H^1({\cal C} ,{\cal C}^{(n)})=0$. Moreover, recall from
\cite{gro} that $\cal C$ is cyclic amenable if $H_\lambda^1({\cal
C} ,{\cal C}^{(1)})=0$. For more information about these notions
we refer the reader to the impressive references \cite{dghg,gro}.


\section{\large\bf Definition and some basic results}

We commence this section by introducing the main object of study of this work.
To this end, we need to recall some notations. Let  ${\mathcal A}$ and ${\frak A}$ be Banach algebras such that ${\mathcal A}$
is a Banach ${\frak A}$-bimodule, we say that the left and right actions of $\frak A$ on
$\mathcal A$ are compatible, if for each $a,b\in{\cal A}$ and $\beta\in{\frak A}$,
\begin{eqnarray*}
\beta\cdot(ab)=(\beta\cdot a)b,\quad\quad(ab)\cdot\beta=a(b\cdot \beta),\quad\quad
a(\beta\cdot b)=(a\cdot\beta)b.
\end{eqnarray*}
and in the case where $\beta\cdot a=a\cdot \beta$, we say that $\cal A$
is a symmetric $\frak A$-bimodule.

\begin{definition}\label{defmain}
Let ${\cal A}$ and ${\frak A}$ be Banach algebras such that ${\cal A}$ is a Banach ${\frak A}$-bimodule with compatible actions. The amalgamated
duplication of $\frak A$ along $\mathcal A$, denoted by ${\cal A}\rtimes{\frak A}$, is defined as the Cartesian product
${\cal A}\times{\frak A}$ with the algebra product
$$(a, \beta)(b, \gamma)={\Big(}ab+a\cdot\gamma+\beta\cdot b,
\beta\gamma{\Big)}$$
and with the norm $\|(a,\beta)\|=\|a\|+\|\beta\|$.
\end{definition}

The following remark is now immediate:

\begin{remark}
{\rm Let ${\cal A}$ and ${\frak A}$ be Banach algebras such that ${\cal A}$ is a Banach ${\frak A}$-bimodule with compatible actions. Then the following statements hold.
\newcounter{j1214}
\begin{list}%
{\rm(\roman{j1214})}{\usecounter{j1214}}
\item It is clear that ${\cal A}\rtimes{\frak A}$ is a Banach algebra. Moreover,
if we identify $\cal A$ with ${\cal A}\times \{0\}$ and $\frak A$ with $\{0\}\times{\frak A}$, then $\cal A$ is a closed ideal
of ${\cal A}\rtimes{\frak A}$ and $\frak A$ is a closed subalgebra.
Moreover
$${\cal A}\rtimes{\frak A}/{\cal A}\cong{\frak A}\quad\quad\quad\quad {\hbox{(isometric~ isomorphism)}}.$$

 This allows us to consider ${\cal A}\rtimes{\frak A}$ as a strongly splitting Banach algebra extension of
$\frak A$ by $\cal A$.

\item ${\cal A}\rtimes{\frak A}$ is commutative if and only if
$\cal A$ and $\frak A$ are commutative Banach algebra and
$\cal A$ is a symmetric $\frak A$-bimodule.

\item If $\mathcal A$ and $\frak A$ are only commutative rings
without Banach space structure, then ${\cal A}\rtimes{\frak A}$
coincides with the amalgamated duplication of a ring $\frak A$
along a $\frak A$-module $\mathcal A$, see \cite[Page 3]{ring2}.
\end{list}}
\end{remark}

The following example shows that the amalgamated duplication of
Banach algebras includes a large class of well-known Banach
algebras.

\begin{example}\label{exmain}
{\rm Let ${\frak A}$ be a Banach algebras. Observe that:
\newcounter{j1215}
\begin{list}%
{\rm(\roman{j1215})}{\usecounter{j1215}} \item If ${\cal A}$ ia a
nonzero Banach ${\frak A}$-bimodule. Then we can consider ${\cal
A}$ as a Banach ${\frak A}$-bimodule with compatible actions when
$\cal A$ is endowed with the zero algebra product. Hence, ${\cal
A}\rtimes {\frak A}$ is the module extension Banach algebra in the
sense of \cite{zhang}. Moreover, as has been discussed in
\cite{zhang}, every triangular Banach algebra is isometrically
isomorphic to a module extension Banach algebra. This allows us to
consider triangular Banach algebras as an example of the
amalgamated duplication of $\frak A$ along $\mathcal A$.

\item If ${\cal A}$ is a Banach algebras, $\sigma({\frak
A})\neq\emptyset$ and $\theta\in\sigma({\frak A})$. Then it is
trivial that with the following actions
$$
\beta\cdot a=a\cdot\beta=\theta(\beta)a\quad\quad\quad
\quad(a\in{\cal A},\beta\in{\frak A}),
$$
${\cal A}$ is a Banach ${\frak A}$-bimodule with compatible actions. Hence, in this case ${\cal A}\rtimes {\frak A}$ is the $\theta$-Lau product ${\mathcal A}\times_\theta{\frak A}$ in the sense of \cite{M}.
\end{list}
}
\end{example}


{\bf Convention}. For the rest of this paper we assume that
${\cal A}$ and ${\frak A}$ are Banach algebras such that ${\cal A}$ is a Banach ${\frak A}$-bimodule with compatible actions, unless stated otherwise.\\


Let $N$ be a left ideal of ${\cal A}\rtimes{\frak A}$ and
$$I_{N}:={\Big\{}a\in{\cal A}:~(a,\beta)\in N~{\hbox{for some}}~\beta\in{\frak A}{\Big\}},$$
$$J_{N}:={\Big\{}\beta\in{\frak A}:~(a,\beta)\in N~{\hbox{for some}}~a\in{\cal A}{\Big\}}.$$
Then a routine computation shows that $J_{N}$ is a left ideal of
$\frak A$ and it is shown by examples (see \cite[Example 2.8]{M})
that in general we can neither expect $I_{N}$ to be an ideal, nor
to have $N=I_{N}\times J_{N}$.

The following result which is stated for left ideals is also true
for right and two-sided ideals. Part (i) of this proposition is a
generalization of two results proved by Monfared
\cite[Propositions 2.6 and 2.7]{M}, and in particular, parts
(ii)-(vi) of this proposition can be considered as new results for
Lau products of Banach algebras. Moreover, as far as we know the
subject, this gives a characterization of the ideal structure of
module extensions of Banach algebras and triangular Banach
algebras.

\begin{proposition}
If $I$ and $J$ are nonempty subspaces of $\cal A$ and $\frak A$,
respectively, then the following assertions hold.
\newcounter{j1216}
\begin{list}%
{\rm(\roman{j1216})}{\usecounter{j1216}} \item The subspace
$I\times J$ is a left ideal of ${\cal A}\rtimes{\frak A}$ if and
only if $I$ and $J$ are left ideals of $\cal A$ and $\frak A$,
respectively, and $I$ is a left $\frak A$-submodule of $\cal A$
for which ${\cal A}\cdot J\subseteq I$.

\item Let $N$ be a left ideal of ${\cal A}\rtimes{\frak A}$
containing $\{0\}\times{\frak A}$. Then $I_{N}$ is a left ideal of $\cal A$ and $N=I_N\times{\frak A}$.

\item Let $N$ be a left ideal of ${\cal A}\rtimes{\frak A}$
containing ${\cal A}\times\{0\}$. Then $J_{N}$ is a left ideal of $\frak A$ and $N={\cal A}\times J_N$.

\item The subspace $I$ is a maximal left ideal and a left $\frak
A$-submodule of $\cal A$ for which ${\cal A}\cdot {\frak
A}\subseteq I$  if and only if $I\times{\frak A}$ is a maximal
left ideal of ${\cal A}\rtimes{\frak A}$.

\item The subspace $J$ is a maximal left ideal of $\frak A$ if and
only if ${\cal A}\times J$ is a maximal left ideal of ${\cal
A}\rtimes{\frak A}$.

\item The subspace $I\times J$ is a maximal left ideal of ${\cal
A}\rtimes{\frak A}$ if and only if either $I={\cal A}$ with $J$ is
a maximal left ideal of ${\frak A}$ or $J={\frak A}$ with $I$ is
both a maximal left ideal and a left $\frak A$-submodule of $\cal
A$ satisfying ${\cal A}\cdot {\frak A}\subseteq I$.
\end{list}
\end{proposition}
{\noindent Proof.} For briefness we only give the proof for assertions (ii), (iv) and (vi).

(ii) Let $N$ be a left ideal of ${\cal A}\rtimes{\frak A}$
containing $\{0\}\times{\frak A}$ and let $a\in I_N$. Then there
is a $\beta$ in ${\frak A}$ such that $(a,\beta)\in N$. It follows
that $(a,0)\in N$. From this, we can deduced that $I_N$ is a left
ideal of $\cal A$ and $N=I_N\times{\frak A}$.

(iv) Let $I$ be a left $\frak A$-submodule and a maximal left ideal of $\cal A$ for which
${\cal A}\cdot {\frak A}\subseteq I$.
That $I\times {\frak A}$ is a left
ideal of ${\cal A}\rtimes {\frak A}$ follows from part (i). Now, let $N$
be any left ideal of ${\cal A}\rtimes {\frak A}$ with $$I\times {\frak A}\subsetneqq N \subseteq {\cal A}\rtimes {\frak A}.$$ Then, by
part (ii),  $I_N$
is a left ideal of $\cal A$ and $N=I_N\times{\frak A}$. It follows that
$I\subsetneqq I_N\subseteq {\cal A}$ and thus $I_N={\cal A}$. Hence $I\times {\frak A}$ is maximal.

Conversely, let $I\times {\frak A}$ be a maximal left ideal of
${\cal A}\rtimes {\frak A}$ and let $I'$
be any left ideal of $\cal A$ for which $I\subsetneqq I'
\subseteq {\cal A}$. Then $$I\times{\frak A}\subsetneqq I'\times{\frak A}
\subseteq {\cal A}\rtimes{\frak A}.$$ As $I\times {\frak A}$ is
maximal, $I'={\cal A}$. Therefore $I$ is a maximal left ideal of $\frak A$.

(vi) Let $I\times J$ be a maximal left ideal of ${\cal A}\rtimes {\frak A}$ and $I\subsetneqq {\cal A}$. Then $$I\times J\subsetneqq {\cal A}\times J\subseteq {\cal A}\rtimes {\frak A}.$$ Since $I\times J$ is a maximal left ideal of ${\cal A}\rtimes {\frak A}$,
we have $J={\frak A}$. Similarly, if $J\subsetneqq {\frak A}$, then
$I={\cal A}$. The converse follows from parts (iv) and (v).$\e$\\


Let ${\mathcal C}$ be a Banach algebra and let  $\cal X$ and $\cal
Y$ be two Banach ${\mathcal C}$-bimodules. An operator $T:{\cal
X}\rightarrow {\cal Y}$ is called a {\it left ${\mathcal
C}$-module map} if $T(c\cdot x) = c\cdot T(x)$ for all
$c\in{\mathcal C}$ and $x\in {\cal X}$. {\it Right ${\mathcal
C}$-module} and {\it ${\mathcal C}$-bimodule maps}  are defined
similarly. We denote by ${\rm Hom}_{\mathcal C}({\cal X},{\cal
Y})$ the space of all bounded left ${\mathcal C}$-module maps from
$\cal X$ into $\cal Y$.  Define ${\rm LM}({\mathcal C}):={\rm
Hom}_ {\mathcal C}({\mathcal C},  {\mathcal C})$ to be the left
multiplier algebra of ${\mathcal C}$.

The following result, stated for left multipliers, is also true
for right multipliers. In particular, it can be considered as a
new result for certain examples of ${\cal A}\rtimes{\frak A}$
which were introduced in Example \ref{exmain}.


\begin{proposition}\label{12iii}
The operator $T$ is in ${\rm LM}({\cal A}\rtimes{\frak A})$ if and only if there exists some
$T^{1}_{\mathcal A}\in {\rm Hom}_{\frak A}({\mathcal A}, {\mathcal A})$,
$T^{1}_{\frak A}\in {\rm Hom}_{\frak A}({\frak A}, {\mathcal A})$,
$T^{2}_{\frak A}\in {\rm LM}({\frak A})$ and $T^{2}_{\mathcal A}\in
{\rm Hom}_{\frak A}({\mathcal A}, {\frak A})$ such that for each
$a, b\in {\mathcal A}$ and $\beta\in {\frak A}$
\newcounter{j1217}
\begin{list}%
{\rm(\roman{j1217})}{\usecounter{j1217}}
\item $T((a,\beta)) = (T^{1}_{\mathcal A}(a) +T^{1}_{\frak A}(\beta),
T^{2}_{\mathcal A}(a) +T^{2}_{\frak A}(\beta))$.

\item $T^{1}_{\mathcal A}(ab)=aT^{1}_{\mathcal A}(b)+a\cdot T^{2}_{\mathcal A}(b)$.

\item $T^{1}_{\mathcal A}(a\cdot\beta)=aT^{1}_{\frak A}
(\beta)+a\cdot T^{2}_{\frak A}(\beta)$.

\item $T^{2}_{\mathcal A}(ab)=T^{2}_{\mathcal A}(a\cdot\beta)=0$.
\end{list}
\end{proposition}
{\noindent Proof.} We only need to prove the ``if" part of this
proposition, which is the essential part of it. Suppose that $T\in
{\rm LM}({\cal A}\rtimes{\frak A})$. Then,
 there exists bounded linear maps $T_1:{\cal A}\rtimes{\frak A}\rightarrow {\mathcal A}$  and $T_2:{\cal A}\rtimes{\frak A}\rightarrow
{\frak A}$  such that $T = (T_1, T_2)$. Let
$$T_{\mathcal A}^1(a)=T_1((a,0)),\quad\quad T_{\mathcal A}^2(a)=T_2((a,0))\quad\quad\quad\quad\quad\quad(a\in{\cal A}),$$ and $$T_{\frak A}^1(\beta)=T_1((0,\beta)),
\quad\quad T_{\frak A}^2(\beta)=T_2((0,\beta))
\quad\quad\quad\quad\quad\quad(\beta\in{\frak A}).$$
Then trivially $T_{\mathcal A}^1$, $T_{\frak A}^1$, $T_{\mathcal A}^2$
and $T_{\frak A}^2$ are linear maps satisfying (i). Moreover, for each
$a,b\in {\mathcal A}$ and $\beta, \gamma\in {\frak A}$, we observe
\begin{eqnarray}\label{1}
T{\Big(}(a,\beta)(b,\gamma){\Big)}&=&T{\Big(}(ab+a\cdot\gamma+
\beta\cdot b,\beta\gamma){\Big)}\nonumber\\
&=&{\Big(}T_{\mathcal A}^1(ab+a\cdot\gamma+\beta\cdot b)+ T_{\frak
A}^1(\beta\gamma),T_{\mathcal A}^2(ab+a\cdot\gamma+\beta\cdot b)
+T_{\frak A}^2(\beta\gamma){\Big)},
\end{eqnarray}
and
\begin{eqnarray}\label{2}
(a,\beta)T((b,\gamma))
&=&(a,\beta){\Big(}T^{1}_{\mathcal A}(b) +T^{1}_{\frak A}(\gamma),
T^{2}_{\mathcal A}(b) +T^{2}_{\frak A}(\gamma){\Big )}\nonumber \\
&=&{\Big(} aT^{1}_{\mathcal A}(b)+aT^{1}_{\frak A}(\gamma)
+a\cdot T^{2}_{\mathcal A}(b)+a\cdot T^{2}_{\frak A}(\gamma)\nonumber\\
&&\quad\quad\quad+\beta\cdot T^{1}_{\mathcal A}(b)+ \beta\cdot
T^{1}_{\frak A}(\gamma),\beta T^{2}_{\mathcal A}(b) +\beta
T^{2}_{\frak A}(\gamma){\Big)}.
\end{eqnarray}
Therefore, Equalities (\ref{1}) and (\ref{2}) imply that
\begin{eqnarray}\label{3}
T_{\mathcal A}^1(ab+a\cdot\gamma+\beta\cdot b)+
T_{\frak A}^1(\beta\gamma)&=&aT^{1}_{\mathcal A}(b)+
aT^{1}_{\frak A}(\gamma)
+a\cdot T^{2}_{\mathcal A}(b)\nonumber\\
&&\quad\quad+a\cdot T^{2}_{\frak A}(\gamma)+\beta\cdot T^{1}_{\mathcal A}(b)+
\beta\cdot T^{1}_{\frak A}(\gamma),
\end{eqnarray}
and
\begin{eqnarray}\label{4}
T_{\mathcal A}^2(ab+a\cdot\gamma+\beta\cdot b)+T_{\frak A}^2(\beta\gamma)=\beta T^{2}_{\mathcal A}(b)+\beta T^{2}_{\frak A}(\gamma).
\end{eqnarray}
Applying (\ref{3}) and (\ref{4}) for suitable values of $a,b,
\beta, \gamma$ show that $T^{1}_{\mathcal A}\in {\rm Hom}_{\frak
A}({\mathcal A}, {\mathcal A})$, $T^{1}_{\frak A}\in {\rm
Hom}_{\frak A}({\frak A}, {\mathcal A})$, $T^{2}_{\frak A}\in {\rm
LM}({\frak A})$ and $T^{2}_{\mathcal A}\in {\rm Hom}_{\frak
A}({\mathcal A}, {\frak A})$ and the Equalities (ii)-(iv) are also
satisfied.$\e$\\


From now on, for a set ${Y}$ in a Banach space,
$\langle {Y} \rangle$ denotes the closed linear
span of ${Y}$ in the space.

\begin{corollary}\label{bi2}
Suppose that either
$\langle{\mathcal A}^2\rangle={\cal A}$ or $\langle{\mathcal A}\cdot{\frak A}\rangle={\cal A}$. Then
$T\in {\rm LM}({\cal A}\rtimes{\frak A})$ if and only if there exists some
$T^{1}_{\mathcal A}\in {\rm Hom}_{\frak A}({\mathcal A}, {\mathcal A})\cap{\rm LM}({\cal A})$,
$T^{1}_{\frak A}\in {\rm Hom}_{\frak A}({\frak A}, {\mathcal A})$ and
$T^{2}_{\frak A}\in {\rm LM}({\frak A})$ such that
$T^{1}_{\mathcal A}(a\cdot\beta)=aT^{1}_{\frak A}
(\beta)+a\cdot T^{2}_{\frak A}(\beta)$ and
$$T((a,\beta)) = (T^{1}_{\mathcal A}(a) +T^{1}_{\frak A}(\beta),
T^{2}_{\frak A}(\beta))$$
for all
$a\in {\mathcal A}$ and $\beta\in {\frak A}$.
\end{corollary}


The first dual space of ${\cal A}\rtimes{\frak A}$ can be identified with
${\cal A}^{(1)}\times {\frak A}^{(1)}$ in the natural way
$${\big<} (a^{(1)},\beta^{(1)}),(a,\beta) {\big>}=\langle a^{(1)},
 a\rangle+\langle \beta^{(1)},\beta\rangle,$$
where $a\in {\cal A},~\beta\in {\frak A},~a^{(1)}\in{\cal
A}^{(1)}$ and $\beta^{(1)}\in {\frak A}^{(1)}$. The dual norm on ${\cal
A}^{(1)}\times{\frak A}^{(1)}$ is of course the maximum norm
$\|(a^{(1)},\beta^{(1)})\|=\max{\Big\{}\|a^{(1)}\|, \|\beta^{(1)}\|{\Big\}}$.

The proof of Proposition \ref{ib} below, which characterize the
spectrum of the Banach algebra ${\cal A}\rtimes{\frak A}$,
relies on the following lemma.

\begin{lemma}\label{231094}
If $\varphi\in \sigma({\mathcal A})$, then there exists a
unique linear functional $\widetilde{\varphi}$ in $\sigma({\frak A})\cup\{0\}$ such that
\begin{eqnarray}\label{5}
\varphi(a\cdot\beta)=\varphi(\beta\cdot a)=\varphi(a)\widetilde{\varphi}(\beta)\quad\quad\quad\quad\quad\quad(a\in{\mathcal A}~ {\hbox{and}}~ \beta\in{\frak A}).
\end{eqnarray}
In particular, if either
$\langle{\mathcal A}\cdot{\frak A}\rangle={\cal A}$ or $\langle{\frak A}\cdot{\mathcal A}\rangle={\cal A}$, then $\widetilde{\varphi}\neq 0$.
\end{lemma}
{\noindent Proof.} Let $\varphi\in\sigma({\mathcal A})$ and $a_0,
a_0'\in{\mathcal A}$ be such that $\varphi(a_0)=\varphi(a_0')=1$.
The compatibility of the left and right actions of $\frak A$ on
$\mathcal A$ imply that for each $\beta\in{\frak A}$,
$$\varphi(\beta\cdot a_0)=\varphi(a_0')\varphi(\beta\cdot a_0)
=\varphi((a_0'\cdot\beta)a_0)=\varphi(a_0'\cdot\beta)
=\varphi(a_0'(\beta\cdot a_0'))
=\varphi(\beta\cdot a_0').$$
We now, define
$\widetilde{\varphi}:{\frak A}\rightarrow{\Bbb C}$ by
$$\widetilde{\varphi}(\beta):=\varphi(\beta\cdot a_0)\quad\quad\quad\quad(\beta\in{\frak A}).$$
The preceding equation shows that the definition is independent of
$a_0$ and hence $\widetilde{\varphi}$ is a well-defined bounded
linear functional on $\frak A$. In particular, we have
$$\varphi(a_0\cdot\beta)=\varphi((a_0\cdot\beta)a_0)=\varphi(\beta\cdot a_0)
\quad\quad\quad(\beta\in{\frak A}).$$ It follows that
$$
\varphi(\beta\cdot a)=\varphi(a_0\cdot\beta)\varphi(a)
=\varphi(\beta\cdot a_0)\varphi(a)=\widetilde{\varphi}(\beta)\varphi(a),
$$
and
$$\varphi(a\cdot\beta)=\varphi(a\cdot(\beta\cdot a_0))=
\varphi(a)\varphi(\beta\cdot
a_0)=\widetilde{\varphi}(\beta)\varphi(a),$$ where $a\in{\mathcal
A}$ and $\beta\in{\frak A}$. Finally, it should be noted that the
uniqueness and multiplicativity of
$\widetilde{\varphi}$ follow similarly.$\e$\\

{\bf Note}:~ {\rm (i) For the rest of this paper, we shall use the letters $\widetilde{\varphi}$
exclusively to denote the unique multiplicative linear functional
associated to $\varphi\in\sigma({\mathcal A})$ which satisfies in (\ref{5}).

(ii) If $\sigma({\mathcal A})\neq\emptyset$ and either
$\langle{\mathcal A}\cdot{\frak A}\rangle={\cal A}$ or $\langle{\frak A}\cdot{\mathcal A}\rangle={\cal A}$, then $\sigma({\frak A})\neq\emptyset$.
}


\begin{proposition}\label{ib}
Let
\begin{eqnarray*}
E:={\Big\{}(\varphi,\widetilde{\varphi}):~ \varphi\in \sigma({\mathcal
A}){\Big\}}\quad{\rm and}\quad F:={\Big\{}(0,\psi):~ \psi\in \sigma({\frak A}){\Big\}}.
\end{eqnarray*}
Set $E=\emptyset$ {\rm[}respectively, $F=\emptyset${\rm]} if
$\sigma({\mathcal A})=\emptyset$ {\rm[}respectively,
$\sigma({\frak A})=\emptyset${\rm]}. Then $E$ and $F$ are disjoint
and $\sigma({\cal A}\rtimes{\frak A})=E\cup F$.
\end{proposition}
{\noindent Proof.} It is clear that $E\cap F=\emptyset$ and $E\cup
F\subseteq \sigma({\cal A}\rtimes{\frak A})$. In order to prove
that $\sigma({\cal A}\rtimes{\frak A})\subseteq E\cup F$, suppose
that $(\varphi,\psi)\in \sigma({\cal A}\rtimes{\frak A})$ and
$(a,\beta), (b,\gamma)\in {\cal A}\rtimes{\frak A}$. Observe that
$$
{\big<}(\varphi,\psi),(a,\beta)(b,\gamma){\big>}={\big<}(\varphi,\psi)
,(a,\beta){\big>}{\big<}(\varphi,\psi),(b,\gamma){\big>}.
$$
It follows that
$$\varphi(ab+a\cdot\gamma+\beta\cdot b)+\psi(\beta\gamma)=
\varphi(a)\varphi(b)+\varphi(a)
\psi(\gamma)+\psi(\beta)\varphi(b)+\psi(\beta)\psi(\gamma).$$
Taking $\beta=\gamma=0$, we have
$\varphi(ab)=\varphi(a)\varphi(b)$ and thus $\varphi\in
\sigma({\mathcal A})\cup\{0\}$. Similarly, we can see that
$\psi\in \sigma({\frak A})\cup\{0\}$. Now, if $\beta=0$ and $b=0$,
then $\varphi(a\cdot\gamma)=\varphi(a)\psi(\gamma)$, similarly
$\varphi(\beta\cdot b)=\psi(\beta)\varphi(b)$.  The equality
$\varphi=0$ implies that $\psi\neq 0$, and if
$\varphi\neq 0$, then $\psi=\widetilde{\varphi}$ by Lemma \ref{231094}.$\e$\\


Next we turn our attention to the question of semisimplicity and
regularity of ${\cal A}\rtimes{\frak A}$.

\begin{proposition}
The following statements hold.
\newcounter{j1218}
\begin{list}%
{\rm(\roman{j1218})}{\usecounter{j1218}}
\item ${\cal A}\rtimes{\frak A}$ is semisimple if and only if
both $\cal A$ and $\frak A$ are semisimple.

\item ${\cal A}\rtimes{\frak A}$ is regular if and only if
both $\cal A$ and $\frak A$ are regular.
\end{list}
\end{proposition}
{\noindent Proof.} (i) Let $\cal A$ and $\frak A$ be semisimple.
We show that the Gelfand homomorphism $\Gamma_{{\cal
A}\rtimes{\frak A}}$ is injective. To this end, suppose that
$a\in{\cal A}$  and $\beta\in{\frak A}$ such that
$\widehat{(a,\beta)}=0$. From this, by Proposition \ref{ib} and
the equality $\widehat{(a,\beta)}= (\hat{a},\hat{\beta})$, we
deduced that
\begin{eqnarray*}
&&0=\big<(\hat{a},\hat{\beta}),(\varphi,\widetilde{\varphi})\big>=
\varphi(a)+\widetilde{\varphi}(\beta)\quad\quad\quad({\varphi\in\sigma({\cal A}))},\\
&&0=\big<(\hat{a},\hat{\beta}),(0,\psi)\big>=
\psi(\beta)\quad\quad\quad\quad\quad\quad~{(}\psi\in\sigma({\frak A}){)}.
\end{eqnarray*}
This implies that $a=\beta=0$. The converse can be proved similarly.

(ii) Since ${\cal A}\times\{0\}$ is a closed
ideal of ${\cal A}\rtimes{\frak A}$ and ${\cal A}\rtimes{\frak A}/{\cal A}\times\{0\}$ is isometrically isomorphic to $\frak A$, it follows from \cite[Theorems 4.2.6 and
4.3.8]{kani} that ${\cal A}\rtimes{\frak A}$ is regular if and only if $\cal A$
and ${\frak A}$ are regular.$\e$\\


Now, let $a, b\in{\cal A}$, $\beta,\gamma\in{\frak A}$, $a^{(1)}\in{\cal A}^{(1)}$,
$a^{(2)}\in{\cal A}^{(2)}$ and $\beta^{(2)}\in{\frak A}^{(2)}$. Then we can extend the left and right actions of $\frak A$ on $\mathcal A$ to compatible actions
of ${\frak A}^{(2)}$ on ${\cal A}^{(2)}$ via
\begin{eqnarray*}
\big<a^{(2)}\bullet\beta^{(2)},a^{(1)}\big>&=&
\big<a^{(2)},\beta^{(2)}\bullet a^{(1)}\big>,\\
\quad\big<\beta^{(2)}\bullet a^{(1)},a\big>&=&
\big<\beta^{(2)},a^{(1)}\bullet a\big>,\\
\big<a^{(1)}\bullet a,\gamma\big>&=&\big<a^{(1)}, a\cdot\gamma\big>,
\end{eqnarray*}
and
\begin{eqnarray*}
\big<\beta^{(2)}\bullet a^{(2)},a^{(1)}\big>&=&
\big<\beta^{(2)},a^{(2)}\bullet a^{(1)}\big>,\\
\big<a^{(2)}\bullet a^{(1)},\beta\big>&=&
\big<a^{(2)},a^{(1)}\bullet\beta\big>,\\
\big<a^{(1)}\bullet\beta,b\big>&=&\big<a^{(1)}, \beta\cdot b\big>.
\end{eqnarray*}
Similarly, by using symmetry, we may consider ${\cal A}^{(2)}$ as
a Banach ${\frak A}^{(2)}$-bimodule with compatible actions via
the following module actions
\begin{eqnarray*}
\big<a^{(2)}\blacktriangle\beta^{(2)},a^{(1)}\big>&=&
\big<\beta^{(2)},a^{(1)}\blacktriangle a^{(2)}\big>,\\
\big<a^{(1)}\blacktriangle a^{(2)},\beta\big>&=&
\big<a^{(2)},\beta\blacktriangle a^{(1)}\big>,\\
\big<\beta\blacktriangle a^{(1)},b\big>&=&
\big<a^{(1)}, b\cdot\beta\big>,
\end{eqnarray*}
and
\begin{eqnarray*}
\big<\beta^{(2)}\blacktriangle a^{(2)},a^{(1)}\big>&=&
\big<a^{(2)},a^{(1)}\blacktriangle\beta^{(2)}\big>,\\
\big<a^{(1)}\blacktriangle \beta^{(2)},a\big>&=&
\big<\beta^{(2)},a\blacktriangle a^{(1)}\big>,\\
\big<a\blacktriangle a^{(1)},\gamma\big>&=&
\big<a^{(1)}, \gamma\cdot a\big>.
\end{eqnarray*}
Hence, one can calculate the first and second Arens products on $({\cal A}\times{\frak A})^{(2)}$ as follows
\begin{eqnarray*}
{\Big(}b^{(1)},\gamma^{(1)}{\Big)}\circ(b,\gamma)&=&{\Big(}b^{(1)}\circ b+b^{(1)}\bullet\gamma,
b^{(1)}\bullet b+\gamma^{(1)}\circ\gamma{\Big)},\\
{\Big(}b^{(2)},\gamma^{(2)}{\Big)}\circ{\Big(}b^{(1)},
\gamma^{(1)}{\Big)}&=&{\Big(}b^{(2)}\circ b^{(1)}
+\gamma^{(2)}\bullet b^{(1)},b^{(2)}\bullet b^{(1)}+\gamma^{(2)}\circ\gamma^{(1)}{\Big)},\\
{\Big(}a^{(2)},\beta^{(2)}{\Big)}\circ{\Big(}b^{(2)},\gamma^{(2)}{\Big)}&=&
{\Big(}a^{(2)}\circ b^{(2)}+a^{(2)}\bullet\gamma^{(2)}+\beta^{(2)}\bullet b^{(2)}
,\beta^{(2)}\circ\gamma^{(2)}{\Big)},
\end{eqnarray*}
and
\begin{eqnarray*}
(b,\gamma)\vartriangle{\Big(}b^{(1)},\gamma^{(1)}{\Big)}&=&{\Big(}b\vartriangle b^{(1)}+\gamma\blacktriangle b^{(1)},b\blacktriangle b^{(1)}+\gamma\vartriangle\gamma^{(1)}{\Big)},\\
{\Big(}b^{(1)},\gamma^{(1)}{\Big)}\vartriangle{\Big(}a^{(2)},\beta^{(2)})&=&
{\Big(}b^{(1)}\vartriangle a^{(2)}
+b^{(1)}\blacktriangle\beta^{(2)},b^{(1)}\blacktriangle a^{(2)}+\gamma^{(1)}\vartriangle\beta^{(2)}{\Big)},\\
{\Big(}a^{(2)},\beta^{(2)}{\Big)}\vartriangle{\Big(}b^{(2)},\gamma^{(2)}{\Big)}&=&
{\Big(}a^{(2)}\vartriangle b^{(2)}+a^{(2)}\blacktriangle\gamma^{(2)}+\beta^{(2)}
\blacktriangle b^{(2)},\beta^{(2)}\vartriangle\gamma^{(2)}{\Big)},
\end{eqnarray*}
where $(b,\gamma)\in{\cal A}\rtimes{\frak A}$,
${\Big(}b^{(1)},\gamma^{(1)}{\Big)}\in ({\cal A}\rtimes{\frak
A})^{(1)}$ and ${\Big(}a^{(2)},\beta^{(2)}{\Big)},
{\Big(}b^{(2)},\gamma^{(2)}{\Big)}\in ({\cal A}\rtimes{\frak
A})^{(2)}$. We summarize these observations in the next result
which shows that the first and second Arens products defined on
$({\cal A}\times{\frak A})^{(2)}$ behaves in a natural way.

\begin{proposition}\label{seconddual}
Suppose that ${\mathcal A}^{(2)}$, ${\frak A}^{(2)}$, and
$({\mathcal A}\rtimes{\frak A})^{(2)}$ are equipped with their
first {\rm(}resp. second{\rm)} Arens products. Then ${\cal
A}^{(2)}$ is a Banach ${\frak A}^{(2)}$-bimodule with compatible
actions and in particular
$$({\mathcal A}\rtimes{\frak A})^{(2)}\cong{\mathcal A}^{(2)}\rtimes{\frak A}^{(2)},$$
where $\cong$ denotes the isometric algebra isomorphism.
\end{proposition}


Let us recall from \cite{esfi} that the topological centres of the left and right
module actions of ${\frak A}^{(2)}$ on ${\cal A}^{(2)}$ are as follows:
$${\frak Z}_{\frak A}({\cal A}^{(2)})={\Big\{}a^{(2)}\in {\cal A}^{(2)}:~
a^{(2)}\bullet\gamma^{(2)}=a^{(2)}\blacktriangle\gamma^{(2)}~~{\hbox{for all}}~~\gamma^{(2)}\in{\frak A}^{(2)}{\Big\}},$$
and
$${\frak Z}_{\cal A}({\frak A}^{(2)})={\Big\{}\beta^{(2)}\in {\frak A}^{(2)}:~
\beta^{(2)}\bullet b^{(2)}=\beta^{(2)}\blacktriangle
b^{(2)}~~{\hbox{for all}}~~b^{(2)}\in{\cal A}^{(2)}{\Big\}}.$$
Moreover, the actions of ${\frak A}$ on ${\cal A}$ are Arens
regular if $${\frak Z}_{\frak A}({\cal A}^{(2)})={\cal A}^{(2)},
\quad\quad{\hbox{and}}\quad\quad{\frak Z}_{\cal A}({\frak
A}^{(2)})={\frak A}^{(2)}$$ and they are strongly Arens irregular
if
$${\frak Z}_{\frak A}({\cal A}^{(2)})={\cal A},
\quad\quad{\hbox{and}}\quad\quad{\frak Z}_{\cal A}({\frak A}^{(2)})={\frak A}.$$

The following proposition is an immediate consequence of
Proposition \ref{seconddual}.

\begin{proposition}
The following equality holds for the first
topological centre of
${\cal A}\rtimes{\frak A}$:
$${\frak Z}_t{\Big(}({\mathcal A}\rtimes{\frak A})^{(2)}{\Big)}=
{\Big(}{\frak Z}_t({\cal A}^{(2)})\cap{\frak Z}_{\frak A} ({\cal
A}^{(2)}){\Big)}\rtimes{\Big(}{\frak Z}_t({\frak A}^{(2)})
\cap{\frak Z}_{\frak A}({\cal A}^{(2)}){\Big)}.$$ In particular,
$({\mathcal A}\rtimes{\frak A})^{(2)}$ is Arens regular
{\rm[}respectively, strongly Arens irregular{\rm]} if and only if
both $\cal A$ and $\frak A$ are Arens regular {\rm[}respectively,
strongly Arens irregular{\rm]} and $\frak A$ act regularly
{\rm(}resp. strongly irregular{\rm)} on $\cal A$.
\end{proposition}


We end this section by noting that, from induction and the previous
analysis of the first Arens product, one
can show that the $({\mathcal A}\rtimes{\frak A})$-bimodule actions
on $({\mathcal A}\rtimes{\frak A})^{(n)}$ are formulated as
follows:
\[ (a,\beta)\vartriangle{\Big(}a^{(n)},\beta^{(n)}{\Big)}=\left\{%
\begin{array}{ccc}
\hspace{-2mm}\vspace{5mm}{\Big(}a\vartriangle a^{(n)}+\beta\blacktriangle a^{(n)}
+a\blacktriangle\beta^{(n)}, \beta\vartriangle\beta^{(n)}{\Big)} & ~~~~\hbox{if}~ n ~\hbox{is even}\\
{\Big(}a\vartriangle a^{(n)}+\beta\blacktriangle a^{(n)}, a\blacktriangle a^{(n)}
+ \beta\vartriangle \beta^{(n)}{\Big)} & ~~~~\hbox{if}~ n ~\hbox{is odd},\\
\end{array}%
\right.\] and
\[ {\Big(}a^{(n)}, \beta^{(n)}{\Big)}\circ (a, \beta)=\left\{%
\begin{array}{ccc}
\hspace{-2mm}\vspace{5mm}{\Big(} a^{(n)}\circ a+a^{(n)}\bullet\beta+
\beta^{(n)}\bullet a,  \beta^{(n)}\circ \beta{\Big)} & ~~~~\hbox{if}~ n ~\hbox{is even} \\
{\Big(}a^{(n)}\circ a+a^{(n)}\bullet\beta,
a^{(n)}\bullet a+\beta^{(n)}\circ\beta{\Big)} & ~~~~\hbox{if}~ n ~\hbox{is odd}
\end{array}%
\right.\] where $(a,\beta)\in{\mathcal A}\times{\frak A}$,
${\Big(}a^{(n)},\beta^{(n)}{\Big)}\in{\cal A}^{(n)}\times {\frak
A}^{(n)}$ and $n\in{\Bbb Z}^+$.

\section{\large\bf Certain Derivations on ${\cal A} \rtimes {\frak A}$}

This section studies the sets $Z^1{(}{\cal A} \rtimes {\frak A},({\cal A}
\rtimes {\frak A})^{(n)}{)}$ and $H^1{(}{\cal A} \rtimes {\frak A},({\cal A} \rtimes {\frak A})^{(n)}{)}$ for all $n\in {\Bbb N}\cup\{0\}$;
in particular, the question when the continuous map $D:{\cal A}\rtimes{\frak A}\rightarrow({\cal A}\rtimes{\frak A})^{(1)}$ is a cyclic derivation.

First, we note that an argument similar to the proof of Proposition \ref{12iii}
gives the following result which characterize the set of all derivations from ${\cal A} \rtimes {\frak A}$ into ${\cal A} \rtimes {\frak A}$-bimodule $({\cal A}
\rtimes {\frak A})^{(n)}$ for the case where $n$ is a odd positive integer number.
The details are omitted.

\begin{proposition}\label{oddd}
A bounded linear map $D : {\cal A} \rtimes {\frak A}\rightarrow ({\cal A} \rtimes {\frak A})^{(2n+1)}$ is a derivation
if and only if there exists derivations $D^1_{\cal A}:{\cal A}\rightarrow{\cal A}^{(2n+1)}$, $D^1_{\frak A}: {\frak A} \rightarrow  {\cal A}^{(2n+1)}$,
$D^2_{\frak A}:{\frak A}\rightarrow{\frak A}^{(2n+1)}$
and a bounded
linear map $D^2_{\cal A}:{\cal A}\rightarrow{\frak A}^{(2n+1)}$ such that
\newcounter{j1219}
\begin{list}%
{\rm(\roman{j1219})}{\usecounter{j1219}}
\item $D((a, \beta))=(D^1_{\cal A}(a) +D^1_{\frak A}(\beta),
D^2_{\cal A}(a) + D^2_{\frak A}(\beta))$,

\item $D^1_{\cal A}(\beta\cdot a)=D^1_{\frak A}(\beta)\circ a
+\beta\blacktriangle D^1_{\cal A}(a)$
and $~D^1_{\cal A}(a\cdot\beta )=a\vartriangle D^1_{\frak A}
(\beta)+D^1_{\cal A}(a)\bullet\beta$

\item $D^2_{\cal A}(\beta\cdot a)=D^1_{\frak A}(\beta)\bullet a+
\beta\vartriangle D^2_{\cal A}(a)$ and
$~D^2_{\cal A}(a\cdot\beta)=a\blacktriangle D^1_{\frak A}(\beta)+D^2_{\cal A}(a)\circ\beta$,

\item $D^2_{\cal A}(ab)=D^1_{\cal A}(a)\bullet b+a\blacktriangle D^1_{\cal A}(b)$,
\end{list}
for all $a, b\in {\cal A}$ and $\beta\in{\frak A}$.
\end{proposition}

Now we examine inner derivation from ${\cal A} \rtimes {\frak A}$ into $({\cal A} \rtimes {\frak A})^{(2n+1)}$.

\begin{proposition}\label{odd2}
Let $D\in Z^1({\cal A} \rtimes {\frak A}, ({\cal A} \rtimes {\frak
A})^{(2n+1)})$ for some $n\in{\Bbb N}$, and let $D^1_{\cal A}$,
$D^1_{\frak A}$, $D^2_{\cal A}$ and $D^2_{\frak A}$ associated to
$D$ as in {\rm Proposition \ref{oddd}}. If $a^{(2n+1)}\in{\cal
A}^{(2n+1)}$ and $\beta^{(2n+1)}\in{\frak A}^{(2n+1)}$, then
$D=\verb"ad"_{{(}a^{(2n+1)}, \beta^{(2n+1)}{)}}$  if and only if
 $D^1_{\cal A}=\verb"ad"_{a^{(2n+1)}}$,
$D^2_{\frak A}=\verb"ad"_{\beta^{(2n+1)}}$, $D^1_{\frak A}=\verb"ad"_{a^{(2n+1)}}$ and
$$D^2_{\cal A}(a)=a\blacktriangle a^{(2n+1)}-a^{(2n+1)}\bullet a$$
for all $a\in {\cal A}$.
\end{proposition}
{\noindent Proof.}
For the proof, we need to note only that if
$D=\verb"ad"_{{(}a^{(2n+1)}, \beta^{(2n+1)}{)}}$ for some
$a^{(2n+1)}\in{\cal A}^{(2n+1)}$ and
$\beta^{(2n+1)}\in{\frak A}^{(2n+1)}$, then
\begin{eqnarray*}
(D^1_{\cal A}(a), D^2_{\cal A}(a))&=&D((a, 0))\\
&=&\verb"ad"_{{(}a^{(2n+1)}, \beta^{(2n+1)}{)}}(a, 0)\\
&=&(a,0)\vartriangle{\Big(}a^{(2n+1)}, \beta^{(2n+1)}{\Big)}-
{\Big(}a^{(2n+1)}, \beta^{(2n+1)}{\Big)}\circ(a,0)\\
&=&{\Big(}a\vartriangle a^{(2n+1)}, a\blacktriangle a^{(2n+1)}{\Big)}-
{\Big(}a^{(2n+1)}\circ a,a^{(2n+1)}\bullet a{\Big)},
\end{eqnarray*}
for all $a\in {\cal A}$. It follows that
$$D^1_{\cal A}(a)=\verb"ad"_{a^{(2n+1)}}(a)\quad\quad {\hbox{and}}\quad\quad
D^2_{\cal A}(a)=a\blacktriangle a^{(2n+1)}-
a^{(2n+1)}\bullet a\quad\quad\quad\quad\quad {(}a\in {\cal A}{)}.$$
Similarly, $$D^2_{\frak A}(\beta)=\verb"ad"_{\beta^{(2n+1)}}(b)\quad\quad{\hbox{ and}}
\quad\quad D^1_{\frak A}(\beta)=\beta\blacktriangle
a^{(2n+1)}-a^{(2n+1)}\bullet \beta\quad\quad\quad\quad\quad(\beta\in{\frak A}),$$
and this completes the proof.$\e$\\


The following corollary is an immediate consequence of Propositions \ref{oddd} and \ref{odd2} above.

\begin{corollary}
Let $n\in{\Bbb N}$ and $T\in {\rm Hom}_{\frak A}({\cal A},{\frak
A}^{(2n+1)})$ be such that $T|_{{\cal A}^2}=0$. Let also $D_T$
from ${\cal A}\rtimes {\frak A}$ into $({\cal A}\rtimes {\frak
A})^{(2n+1)}$ be defined by $D_T((a,\beta))=(0,T(a))$. Then $D_T$
is in $Z^1({\cal A} \rtimes {\frak A},({\cal A}\rtimes {\frak
A})^{(2n+1)})$. In particular, $D_T$ is inner if and only if there
exists $a^{(2n+1)}\in {\cal A}^{(2n+1)}$ such that $a\vartriangle
a^{(2n+1)}=a^{(2n+1)}\circ a$, $\beta\blacktriangle
a^{(2n+1)}=a^{(2n+1)}\bullet\beta$ and
$$T(a)=a\blacktriangle a^{(2n+1)}-a^{(2n+1)}\bullet a,$$
for all $a\in{\cal A}$ and $\beta\in{\frak A}$.
\end{corollary}


Also, Proposition \ref{oddd} paves the way for characterizing the
cyclic derivations on  ${\cal A}\rtimes{\frak A}$ as follows:

\begin{proposition}\label{odd}
A continuous linear map $D:{\cal A}\rtimes{\frak A}
\rightarrow({\cal A}\rtimes{\frak A})^{(1)}$
is a cyclic derivation if and only if
there exists cyclic derivations $D^1_{\cal A}\in Z^1({\cal A},
{\cal A}^{(1)})$ and $D^2_{\frak A}\in
Z^1({\frak A}, {\frak A}^{(1)})$, a derivation $D^1_{\frak A}\in Z^1({\frak A},
{\cal A}^{(1)})$ and a bounded linear map $D^2_{\cal A}: {\cal A}
\rightarrow {\frak A}^{(1)}$ such that for each $a, b\in {\cal A}$ and $\beta\in{\frak A}$
\newcounter{j12141}
\begin{list}%
{\rm(\roman{j12141})}{\usecounter{j12141}}
\item $D((a, \beta)) = (D^1_{\cal A}(a) + D^1_{\frak A}(\beta),
D^2_{\cal A}(a) + D^2_{\frak A}(\beta))$,

\item $D^1_{\cal A}(\beta\cdot a)=D^1_{\frak A}(\beta)\circ a+
\beta\blacktriangle D^1_{\cal A}(a)$ and
$~D^1_{\cal A}(a\cdot\beta )=a\vartriangle D^1_{\frak A}(\beta)+ D^1_{\cal A}(a)\bullet\beta$,

\item $D^2_{\cal A}(\beta\cdot a)=D^1_{\frak A}(\beta)\bullet a +
\beta\vartriangle D^2_{\cal A}(a)$ and
$~D^2_{\cal A}(a\cdot\beta) =a\blacktriangle D^1_{\frak A}(\beta)+D^2_{\cal A}(a)\circ\beta$,

\item $D^2_{\cal A}(ab) = D^1_{\cal A}(a)\bullet b+a\blacktriangle D^1_{\cal A}(b)$,

\item $(D^2_{\cal A})^*\circ \Phi+D^1_{\frak A}=0$, where $\Phi$ is the usual
injection from $\frak A$ into ${\frak A}^{(2)}$ and $(D^2_{\cal A})^*$ is the adjoint of the operator
$D^2_{\cal A}$.
\end{list}
\end{proposition}
{\noindent Proof.} First assume that $D$ is a cyclic derivation
from ${\cal A}\rtimes{\frak A}$ into $({\cal A}\rtimes{\frak
A})^{(1)}$. By Proposition \ref{oddd} there exists $D^1_{\cal
A}\in Z^1({\cal A}, {\cal A}^{(1)})$, $D^2_{\frak A}\in Z^1({\frak
A}, {\frak A}^{(1)})$, $D^1_{\frak A}\in Z^1({\frak A}, {\cal
A}^{(1)})$ and a bounded linear map $D^2_{\cal A}:{\cal
A}\rightarrow {\frak A}^{(1)}$ satisfying  the conditions
(i)-(iv). Recall from the proof of Proposition \ref{oddd} that if
$D=(D_1,D_2)$, then
$$D^1_{\cal A}(a)=D_1((a,0)),\quad\quad D^2_{\cal A}(a)=D_2((a,0))\quad\quad\quad\quad\quad(a\in{\cal A}),$$ and $$D^1_{\frak A}(b)=D_1((0,\beta)),
\quad\quad D^2_{\frak A}(\beta)=D_2((0,\beta))\quad\quad\quad\quad\quad(\beta\in{\frak A}).$$
Moreover, for each $(a,\beta), (b,\gamma)$ in ${\cal A}\rtimes{\frak A}$, we have
\begin{eqnarray}\label{10}
\big<D((a,\beta)),(b,\gamma)\big>+\big<D((b,\gamma)),(a,\beta)\big>=0.
\end{eqnarray}
 It
follows that
\begin{eqnarray}\label{11}
\big<D^1_{\cal A}(a),b\big>+\big<D^1_{\cal A}(b),a\big>=
\big<D((a,0)),(b,0)\big>+\big<D((b,0)),(a,0)\big>=0,
\end{eqnarray}
 and
\begin{eqnarray}\label{12}
\big<D^2_{\frak A}(\beta),\gamma\big>+\big<D^2_{\frak A}(\gamma),\beta\big>=
\big<D((0,\beta)),(0,\gamma)\big>+\big<D((0,\gamma)),(0,\beta)\big>=0.
\end{eqnarray}
 Hence
$D^1_{\cal A}$ and $D^2_{\frak A}$ are cyclic derivations. Now, by equalities
(\ref{10})-(\ref{12}), we deduce that
\begin{eqnarray}\label{13}
\big<D^2_{\cal A}(a),\gamma\big>+
\big<D^2_{\cal A}(b),\beta\big>+
\big<D^1_{\frak A}(\beta),b\big>+\big<D^1_{\frak A}(\gamma),a\big>=0.
\end{eqnarray}
Choosing $\beta=0$ in (\ref{13}), we see that
$(D^2_{\cal A})^*\circ \Phi+D^1_{\frak A}=0$. Finally, we note that the proof of the
converse is not difficult and is omitted.
$\e$\\


The following proposition is now immediate:

\begin{proposition}\label{cyclic2}
Let $D\in Z^1({\cal A} \rtimes {\frak A}, ({\cal A} \rtimes {\frak A})^{(1)})$ be cyclic, and let $D^1_{\cal A}$, $D^1_{\frak A}$, $D^2_{\cal A}$ and $D^2_{\frak A}$ associated to $D$ as in {\rm Proposition \ref{odd}}. If
$a^{(1)}\in{\cal A}^{(1)}$ and
$\beta^{(1)}\in{\frak A}^{(1)}$, then
$D=\verb"ad"_{{(}a^{(1)}, \beta^{(1)}{)}}$ if and only if $D^1_{\cal A}=\verb"ad"_{a^{(1)}}$,
$D^2_{\frak A}=\verb"ad"_{\beta^{(1)}}$, $D^1_{\frak A}=\verb"ad"_{a^{(1)}}$  and
$D^2_{\cal A}(a)=a\blacktriangle a^{(1)}-
a^{(1)}\bullet a$
for all $a\in {\cal A}$.
\end{proposition}


By an argument similar to the proof of Propositions \ref{oddd} and
\ref{odd2} one can obtain the following result which gives a
characterization for the set of all derivation from ${\cal
A}\times{\frak A}$ into $({\cal A}\times{\frak A})^{(2n)}$. The
details are omitted.

\begin{proposition}\label{even}
A bounded linear map $D:{\cal A}\rtimes{\frak A}\rightarrow({\cal A}\rtimes{\frak A})^{(2n)}$ is a derivation if and only if there exists derivations
$D^2_{\frak A}:{\frak A} \rightarrow  {\frak A}^{(2n)}$, $D^1_{\frak A} : {\frak A} \rightarrow  {\cal A}^{(2n)}$, $D^2_{\cal A}\in {B}_{\frak A}({\cal A},{\frak A}^{(2n)})$
and a bounded linear map
$D^1_{\cal A} : {\cal A} \rightarrow  {\cal A}^{(2n)}$  such that
\newcounter{j12144}
\begin{list}%
{\rm(\roman{j12144})}{\usecounter{j12144}}
\item $D((a, \beta))=(D^1_{\cal A}(a) + D^1_{\frak A}(\beta),
D^2_{\cal A}(a) + D^2_{\frak A}(\beta))$,

\item $D^1_{\cal A}(ab)=D^1_{\cal A}(a)\circ b+D^2_{\cal A}(a)\bullet b
+a\vartriangle D^1_{\cal A}(b)+a\blacktriangle D^2_{\cal A}(b)$,

\item $D^1_{\cal A}(\beta\cdot a )=D^1_{\frak A}(\beta)\circ a+
D^2_{\frak A}(\beta)\bullet a+\beta\blacktriangle D^1_{\cal A}(a)$,

\item $D^1_{\cal A}(a\cdot \beta)=a\vartriangle D^1_{\frak A}
(\beta)+a\blacktriangle D^2_{\frak A}(\beta)+D^1_{\cal A}(a)\bullet\beta$,

\item $~D^2_{\cal A}(ab) =0$,
\end{list}
for all $a, b\in {\cal A}$ and $\beta, \gamma\in{\frak A}$.  Moreover, if $a^{(2n)}\in {\mathcal A}^{(2n)}$ and $\beta^{(2n)}\in{\frak A}^{(2n)}$, then $D=\verb"ad"_{(a^{(2n)}, \beta^{(2n)}{)}}$ if and only if $D^1_{\frak A}= \verb"ad"_{a^{(2n)}}$, $D^2_{\frak A}= \verb"ad"_{\beta^{(2n)}}$, $D^2_{\cal A}=0$ and
$$D^1_{\cal A}(a)=\verb"ad"_{a^{(2n)}}(a)+a\blacktriangle\beta^{(2n)}-\beta^{(2n)}\bullet a,$$ for all $a\in {\mathcal A}$.
\end{proposition}


The following two corollaries are immediate consequences of Proposition \ref{even} above.
Recall that, for each $n\in {\Bbb N}\cup\{0\}$, a bounded linear map
$D:{\cal A}\rightarrow{\cal A}^{(2n)}$ is called a $\frak A$-module derivation if
$$D(ab)=D(a)\circ b+a\vartriangle D(b)
\quad\quad\quad\quad{\Big(}a, b\in{\cal A}{\Big)},$$
and
$$D(\beta\cdot a )=\beta\blacktriangle D(a),\quad\quad
D(a\cdot \beta)=D(a)\bullet\beta
\quad\quad\quad\quad~~~~~~~~~~~~~~{\Big(}\beta\in{\frak A}, a\in{\cal A}{\Big)}.$$
The notation $Z^1_{\frak A}({\cal A},{\cal A}^{(2n)})$ is used to denote
the set of all module derivations from ${\cal A}$ into ${\cal A}^{(2n)}$.

\begin{corollary}
Suppose that $T\in Z^1_{\frak A}({\cal A},{\cal A}^{(2n)})$ for some positive integer $n$. Then
the bounded linear map $D_T$ defined by
$$D_T:{\cal A}\rtimes{\frak A}\rightarrow({\cal A}\rtimes{\frak A})^{(2n)};~D_T((a,\beta))=(T(a),0),$$ is in $Z^1({\cal A}\rtimes{\frak A},({\cal A}\rtimes{\frak A})^{(2n)})$. Moreover,  $D_T$ is inner if and only if
there exists $a^{(2n)}\in {\cal A}^{(2n)}$ such that $\beta\blacktriangle a^{(2n)}=a^{(2n)}\bullet \beta$ and
$T(a)=\verb"ad"_{a^{(2n)}}(a)$ for all $a\in {\cal A}$ and $\beta\in{\frak A}$.
\end{corollary}

\begin{corollary}
Suppose that $D_{\cal A}\in Z^1({\cal A},{\cal A}^{(2n)})$ and $T\in Z^1_{\frak A}({\cal A},{\cal A}^{(2n)})$ for some positive integer $n$. Then $D_T:{\cal A}\rtimes{\cal A}\rightarrow({\cal A}\rtimes{\cal A})^{(2n)}$ defined by $$D_T((a,b))=(T(a),D_{\cal A}(b))$$ is in $Z^1({\cal A}\rtimes{\cal A},({\cal A}\rtimes{\cal A})^{(2n)})$. Moreover, $D_T$ is inner if and only if
there exists $a^{(2n)}\in {\cal A}^{(2n)}$ such that $\beta\blacktriangle a^{(2n)}=a^{(2n)}\bullet \beta$ and
$T(a)=\verb"ad"_{a^{(2n)}}(a)$ for all $a\in {\cal A}$ and $\beta\in{\frak A}$.
\end{corollary}


Finally, the results of this section lead to the following corollary which characterize
the sets $Z^1{(}{\cal A} \rtimes {\frak A},({\cal A}
\rtimes {\frak A})^{(n)}{)}$ and $H^1{(}{\cal A} \rtimes {\frak A},({\cal A} \rtimes {\frak A})^{(n)}{)}$ in the case where $n\in {\Bbb N}\cup\{0\}$ and ${\cal A}$ is a unital
Banach algebra.

\begin{corollary}\label{cor1}
Let ${\cal A}$ be unital with the identity $e$. Then the following statements hold.
\newcounter{j12145}
\begin{list}%
{\rm(\roman{j12145})}{\usecounter{j12145}}
\item Every derivation
$D : {\cal A} \rtimes {\frak A}\rightarrow ({\cal A} \rtimes {\frak A})^{(2n+1)}$
is in the form of
 $$D((a, \beta)) = (D^1_{\cal A}(a) + D^1_{\frak A}(\beta), D^2_{\cal A}(a) + D^2_{\frak A}(\beta))\quad\quad\quad(a\in {\mathcal A}, \beta\in{\frak A}),$$ where
  $D^1_{\cal A} : {\cal A} \rightarrow  {\cal A}^{(2n+1)}$,
$D^2_{\frak A} : {\frak A} \rightarrow  {\frak A}^{(2n+1)}$, $D^1_{\frak A} : {\frak A} \rightarrow  {\cal A}^{(2n+1)}$
are derivations  and $D^2_{\cal A} : {\cal A} \rightarrow  {\frak A}^{(2n+1)}$ is a bounded linear map  satisfying
 $$D^1_{\frak A}(\beta)=D^1_{\cal A}(e\cdot\beta)=D^1_{\cal A}(\beta\cdot e)\quad{\hbox {and}}\quad D^2_{\cal A}(a)=D^1_{\cal A}(a)\bullet e=e\blacktriangle D^1_{\cal A}(a).$$

\item Every derivation $D : {\cal A} \rtimes {\frak A}\rightarrow ({\cal A} \rtimes {\frak A})^{(2n)}$
is in the form of
$$D((a, \beta)) = (D^1_{\cal A}(a) + D^1_{\frak A}(\beta),
D^2_{\frak A}(\beta))\quad\quad\quad(a\in {\mathcal A}, \beta\in{\frak A}),$$ where
$D^1_{\cal A} : {\cal A} \rightarrow  {\cal A}^{(2n)}$,
$D^2_{\frak A} : {\frak A} \rightarrow  {\frak A}^{(2n)}$, $D^1_{\frak A}: {\frak A} \rightarrow  {\cal A}^{(2n+1)}$
are derivations  satisfying
 $$D^1_{\frak A}(\beta)=D^1_{\cal A}(e\cdot\beta)-D^2_{\frak A}(\beta)\bullet e=D^1_{\cal A}(\beta\cdot e)-e\blacktriangle D^2_{\frak A}(\beta).\\$$
\end{list}
\end{corollary}


\section{\large\bf n-weak amenability and cyclic amenability}

Our aim in this section is to investigate some notions of
amenability for ${\cal A}\rtimes{\frak A}$ in relation to the
corresponding ones of $\cal A$ and $\frak A$. We commence with the
following result which gives a sufficient condition for
$(2n+1)$-weak amenability of ${\cal A}\rtimes{\frak A}$.

\begin{proposition}\label{ddd}
Suppose that $n$ is a
positive integer number for which $\langle{\cal A}\vartriangle{\cal A}^{(2n)}\rangle={\cal A}^{(2n)}$ or $\langle{\cal A}^{(2n)}\circ {\cal A}\rangle={\cal A}^{(2n)}$.
If $\cal A$ and $\frak A$ are $(2n+1)$-weakly amenable, then
${\cal A}\rtimes{\frak A}$ is $(2n+1)$-weakly amenable.
\end{proposition}
{\noindent Proof.} Suppose that $D\in Z^1({\cal A}\rtimes{\frak A},({\cal A}\rtimes{\frak A})^{(2n+1)})$. Then, there exists $D_{\cal A}^1\in Z^1({\cal A},{\cal A}^{(2n+1)})$,
$D_{\frak A}^1\in Z^1({\frak A},{\cal A}^{(2n+1)})$, $D_{\frak A}^2\in Z^1({\frak A},{\frak A}^{(2n+1)})$ and a bounded linear map $D_{\cal A}^2\in Z^1({\cal A},{\frak A}^{(2n+1)})$
such that
$$D((a, \beta))=(D^1_{\cal A}(a) +D^1_{\frak A}(\beta),
D^2_{\cal A}(a) + D^2_{\frak A}(\beta)),$$
for all $a\in{\cal A}$ and $\beta\in{\frak A}$. By assumption $D_{\cal A}^1=\verb"ad"_{a^{(2n+1)}}$
and $D_{\frak A}^2=\verb"ad"_{\beta^{(2n+1)}}$ for some $a^{(2n+1)}\in{\cal A}^{(2n+1)}$
and $\beta^{(2n+1)}\in{\frak A}^{(2n+1)}$. From this, by Proposition \ref{oddd}, we have
\begin{eqnarray*}
D_{\cal A}^2(ab)&=&\verb"ad"_{a^{(2n+1)}}(a)\bullet b+a\blacktriangle\verb"ad"_{a^{(2n+1)}}(b)\\
&=&(a\vartriangle a^{(2n+1)}-a^{(2n+1)}\circ a)\bullet b+a\blacktriangle(b\vartriangle a^{(2n+1)}-a^{(2n+1)}\circ b)\\
&=&a\blacktriangle(b\vartriangle a^{(2n+1)})-(a^{(2n+1)}\circ a)\bullet b\\
&=&(ab)\blacktriangle a^{(2n+1)})-a^{(2n+1)}\bullet (ab),
\end{eqnarray*}
for all $a, b\in{\cal A}$. Moreover, $(2n+1)$-weak amenability of  $\cal A$ implies  that
$\langle{\cal A}^2\rangle={\cal A}$. Therefore,
$$D_{\cal A}^2(a)=a\blacktriangle a^{(2n+1)}-a^{(2n+1)}\bullet a,$$
On the other hand,  by Propositions \ref{oddd} and \ref{odd2}, we have
$$\verb"ad"_{a^{(2n+1)}}(\beta\cdot a)=D_{\frak A}^1(\beta)\circ a+\beta\blacktriangle
\verb"ad"_{a^{(2n+1)}}(a).$$
It follows that
$${\Big(}D_{\frak A}^1(\beta)-(\beta\blacktriangle a^{(2n+1)}-a^{(2n+1)}\bullet\beta){\Big)}\circ a=0.$$
We now invoke the equality
$\langle{\cal A}\vartriangle{\cal A}^{(2n)}\rangle={\cal A}^{(2n)}$ to conclude that
$$D_{\frak A}^1(\beta)=\beta\blacktriangle a^{(2n+1)}-a^{(2n+1)}\bullet\beta.$$
Similarly, if $\langle{\cal A}^{(2n)}\circ {\cal A}\rangle={\cal A}^{(2n)}$,  using  the equality
$$\verb"ad"_{a^{(2n+1)}}(a\cdot\beta)=a\vartriangle D_{\frak A}^1(\beta)+
\verb"ad"_{a^{(2n+1)}}(a)\bullet\beta,$$ we can show that the
equality $D_{\frak A}^1(\beta)=\beta\blacktriangle
a^{(2n+1)}-a^{(2n+1)}\bullet\beta$ holds. Hence,
$D=\verb"ad"_{(a^{(2n+1)},\beta^{(2n+1)})}$
by Proposition \ref{odd2}.$\e$\\


There does not seem to be an easy to show that if ${\cal A}\rtimes{\frak A}$ is $(2n+1)$-weakly amenable, then $\cal A$ is also. We believe that more information about the properties of $\cal A$, $\frak A$ and the actions of $\frak A$ on $\cal A$ are needed for further conclusions. Hence,
motivated by Proposition \ref{oddd}, it seems valuable to define the following property
for the pair $({\cal A},{\frak A})$.

\begin{definition}
{\rm We say that the pair {\rm(}${\cal A}$,
${\frak A}${\rm)} has the property {\rm(}${\Bbb H}_{2n+1}${\rm)} if for each
$D^1_{\cal A}\in Z^1({\cal A},{\cal A}^{(2n+1)})$, there are $D^1_{\frak A}\in Z^1({\frak A},{\cal A}^{(2n+1)})$
and a bounded
linear map
$D^2_{\cal A}:{\cal A}\rightarrow{\frak A}^{(2n+1)}$  satisfying the  conditions (ii)-(iv) of  Proposition \ref{oddd}.}
\end{definition}

Actually the class of the pair {\rm(}${\cal A}$,${\frak A}${\rm)}
satisfying this property is quite rich. It contains for instance
all the pairs {\rm(}${\cal A}$,${\frak A}${\rm)} such that
\newcounter{j12148}
\begin{list}%
{\rm(\roman{j12148})}{\usecounter{j12148}}
\item ${\cal A}={\frak A}$ and the compatible actions is the natural actions;

\item $\cal A$ is unital;

\item $\cal A$ has a bounded approximate identity and
$\cal A$ is an essential Banach $\frak A$-bimodule;

\item  the  compatible actions of ${\frak A}$ on ${\mathcal A}$ are defined as follows:
$$
\beta\cdot a=a\cdot\beta=\theta(\beta)a\quad(a\in{\cal A},\beta\in{\frak A}),
$$
where $\theta\in\sigma({\frak A})$;

\item the  compatible actions of ${\frak A}$ on ${\mathcal A}$ are defined as follows:
$$
\beta\cdot a=T(\beta)a,\quad a\cdot\beta=aT(\beta)
\quad(a\in{\cal A},\beta\in{\frak A}),
$$
where $T$ is a norm decreasing homomorphism  from ${\frak A}$ into $\cal A$.
\end{list}


\begin{theorem}\label{oddweak}
If ${\cal A}\rtimes{\frak A}$ is $(2n+1)$-weakly amenable for some $n\in {\Bbb N}\cup\{0\}$, then the following assertions hold.
\newcounter{j12156}
\begin{list}%
{\rm(\roman{j12156})}{\usecounter{j12156}}
\item $\frak A$ is $(2n+1)$-weakly amenable.

\item If {\rm(}${\cal A}$,${\frak A}${\rm)} has the property
{\rm(}${\Bbb H}_{2n+1}${\rm)}, then
${\cal A}$ is $(2n+1)$-weakly amenable.
\end{list}
\end{theorem}
{\noindent Proof.} (i) Suppose that $D_{\frak A}^2\in Z^1({\frak A},{\frak A}^{(2n+1)})$. Then, by Proposition \ref{oddd}, the map $D:{\cal A}\rtimes{\frak A}\rightarrow({\cal A}\rtimes{\frak A})^{(2n+1)}$ defined by
$$D((a,\beta))=(0,D_{\frak A}^2(\beta))\quad\quad\quad{\Big(}(a,\beta)\in{\cal A}\rtimes{\frak A}{\Big)}$$
is a  derivation and therefore it is inner by $(2n+1)$-weak amenability of ${\cal A}\rtimes{\frak A}$. We now invoke Proposition \ref{odd2} to conclude that $D^2_{\frak A}$
is inner.

(ii) Suppose that $D_{\cal A}^1\in Z^1({\cal A},{\cal
A}^{(2n+1)})$. Since the pair {\rm(}${\cal A}$,${\frak A}${\rm)}
has the property ${\Bbb H}_{2n+1}$, there exists $D_{\frak A}^1\in
Z^1({\frak A},{\cal A}^{(2n+1)})$ and a bounded linear map
$D^2_{\cal A}:{\cal A}\rightarrow{\frak A}^{(2n+1)}$ satisfying
the conditions (ii)-(iv) of {\rm Proposition \ref{oddd}}. From
this, by Proposition \ref{oddd}, we can conclude that the map
$D:{\cal A}\rtimes{\frak A}\rightarrow({\cal A}\rtimes{\frak
A})^{(2n+1)}$ defined by
$$D((a,\beta))=(D^1_{\cal A}(a) +D^1_{\frak A}(\beta),
D^2_{\cal A}(a))\quad\quad\quad{\Big(}(a,\beta)\in{\cal A}\rtimes{\frak A}{\Big)}$$
is a  derivation and therefore it is inner by $(2n+1)$-weak amenability of ${\cal A}\rtimes{\frak A}$. This together with Proposition \ref{odd2} implies that
$D^1_{\cal A}$ is also inner and this completes the proof.$\e$\\


The next result gives a necessary and sufficient condition for $n$-weak
amenability of ${\cal A}\rtimes{\frak A}$
in relation to that of the $n$-weak amenability of  $\cal A$ and $\frak A$ for the case where $\cal A$ is unital.

\begin{theorem}
Let ${\cal A}$ be unital  and let $n\in {\Bbb N}\cup\{0\}$. Then  ${\cal A} \rtimes {\frak A}$ is $n$-weakly
amenable if and only if ${\cal A}$ and ${\frak A}$ are $n$-weakly
amenable.
\end{theorem}
{\noindent Proof.} Let $D : {\cal A} \rtimes {\frak A}\rightarrow ({\cal A} \rtimes {\frak A})^{(2n+1)}$
be a derivation and  let  $D^1_{\cal A}$, $D^1_{\frak A}$, $D^2_{\cal A}$ and $D^2_{\frak A}$  associated to $D$ as in Corollary \ref{cor1}. Then it is not hard to check  that    if $D^1_{\cal A}={\verb"ad"}_{a^{(2n+1)}}$, then
$D^1_{\frak A}={\verb"ad"}_{a^{(2n+1)}}$ and $$D^2_{\cal A}(a)=a\blacktriangle{a^{(2n+1)}}-{a^{(2n+1)}}\bullet a,$$ for all $a\in{\mathcal A}$. It follows that $D={\verb"ad"}_{(a^{(2n+1)}, \beta^{(2n+1)})}$
if and only if $D^1_{\cal A}={\verb"ad"}_{a^{(2n+1)}}$ and $D^2_{\frak A}={\verb"ad"}_{\beta^{(2n+1)}}$ for some $\beta^{(2n+1)}\in{\frak A}^{(2n+1)}$ and $a^{(2n+1)}\in {\cal A}^{(2n+1)}$. Thus ${\cal A} \rtimes {\frak A}$ is $(2n + 1)$-weakly
amenable if and only if ${\cal A}$ and ${\frak A}$ are $(2n+1)$-weakly amenable. A similar
argument is true for even $n$.$\e$\\


\begin{theorem}\label{cyclic}
If $\cal A$ and $\frak A$
are cyclicly amenable and $\langle{\cal A}^2\rangle={\cal A}$, then
${\cal A}\rtimes{\frak A}$ is
cyclicly amenable.
\end{theorem}
{\noindent Proof.}
Suppose that $D:{\cal A}\rtimes{\frak A}\rightarrow({\cal A}\rtimes{\frak A})^{(1)}$ is a cyclic derivation. Then, there exists $D^1_{\frak A} \in
Z^1({\frak A}, {\cal A}^{(1)})$, a bounded
linear map $D_{\cal A}^2: {\cal A} \rightarrow  {\frak A}^{(1)}$ and
two cyclic derivations $D^1_{\cal A} \in Z^1({\cal A}, {\cal A}^{(1)})$,
$D^2_{\frak A} \in Z^1({\frak A}, {\frak A}^{*})$  such that $(D_{\cal A}^2)^*(\beta)+D^1_{\frak A}(\beta)=0$
and
 $$D((a,\beta)) = (D^1_{\cal A}(a) + D^1_{\frak A}(\beta),D_{\cal A}^2(a) + D^2_{\frak A}(\beta)),$$
for all $a\in{\cal A}$ and $\beta\in {\frak A}$.
Moreover, by assumption there exists
$a^{(1)}\in {\cal A}^{(1)}$ and
$\beta^{(1)}\in {\frak A}^{(1)}$ such that
$D^1_{\cal A}=\verb"ad"_{a^{(1)}}$ and
$D^2_{\frak A}=\verb"ad"_{\beta^{(1)}}$. Hence, for each
$a, b\in{\cal A}$, we have
\begin{eqnarray*}
\big<D^1_{\frak A}(\beta),ab\big>&=&-\big<(D_{\cal A}^2)^*(\beta),ab\big>\\
&=&-\big<D_{\cal A}^2(ab),\beta\big>\\
&=&-\big<D_{\cal A}^1(a)\bullet b+a\blacktriangle D_{\cal A}^1(b),\beta\big>\\
&=&-\big<\verb"ad"_{a^{(1)}}(a)\bullet b+a\blacktriangle \verb"ad"_{a^{(1)}}(b),\beta\big>\\
&=&-\big<a\vartriangle a^{(1)}-a^{(1)}\circ a, b\cdot\beta\big>
-\big<b\vartriangle a^{(1)}-a^{(1)}\circ b,\beta\cdot a\big>\\
&=&-\big<a\vartriangle a^{(1)},b\cdot\beta\big>+\big<a^{(1)}\circ a,b\cdot\beta\big>
-\big<b\vartriangle a^{(1)},\beta\cdot a\big>+
\big<a^{(1)}\circ b,\beta\cdot a\big>\\
&=&\big<a^{(1)}\circ a,b\cdot\beta\big>
-\big<b\vartriangle a^{(1)},\beta\cdot a\big>\\
&=&\big<a^{(1)},(ab)\cdot\beta\big>-\big<a^{(1)},\beta\cdot(ab)\big>\\
&=&\big<\beta\blacktriangle a^{(1)},ab\big>-\big<a^{(1)}\bullet\beta,ab\big>\\
&=&\big<\beta\blacktriangle a^{(1)}-a^{(1)}\bullet\beta,ab\big>.
\end{eqnarray*}
This together with the fact that $\langle{\cal A}^2\rangle={\cal A}$, implies that
$$D_{\frak A}^1(\beta)=\beta\blacktriangle a^{(1)}-a^{(1)}\bullet\beta,$$
for all $\beta\in{\frak A}$.
The same argument shows that
$$D^2_{\cal A}(a)=a\blacktriangle a^{(1)}-a^{(1)}\bullet a,$$
for all $a\in{\cal A}$. We now invoke Proposition \ref{cyclic2}
to conclude that $D=\verb"ad"_{(a^{(1)},\beta^{(1)})}$.
$\e$\\


Recall from \cite[Example 2.5]{gro} that if $0\neq {\cal C}$
is a Banach algebra with zero algebra product, then ${\cal C}$
is cyclicly amenable if and only if $\dim{\cal C}=1$. In
particular, cyclic amenability of the Banach
algebra ${\cal C}$ does not implies that
$\langle{\cal C}^2\rangle={\cal C}$.

Now, we are in position to show that the condition
$\langle{\cal C}^2\rangle={\cal C}$ in Theorem \ref{cyclic}
can not be omitted.

\begin{example}
{\rm Let ${\cal A}={\Bbb C}={\frak A}$ with zero algebra products. Then by
\cite[Example 2.5]{gro} $\cal A$ and $\frak A$ are cyclicly
amenable.  In particular, $\langle{\cal A}^2\rangle\subsetneqq{\cal A}$.
But, again by another application of \cite[Example 2.5]{gro},
${\cal A}\rtimes {\frak A}$ is not  cyclicly amenable. This is
because of, ${\cal A}\rtimes {\frak A}$ is a commutative Banach algebra
with zero algebra product and dimension 2.}
\end{example}


We note that with an argument similar to the proof of Theorem \ref{oddweak} one can prove
the following result. The details are omitted.

\begin{theorem}\label{od2}
If ${\cal A}\rtimes{\frak A}$ is cyclicly amenable, then the following assertions hold.
\newcounter{j12154}
\begin{list}%
{\rm(\roman{j12154})}{\usecounter{j12154}}
\item $\frak A$ is cyclicly amenable.

\item If {\rm(}${\cal A}$,${\frak A}${\rm)} has the property
{\rm(}${\Bbb H}_{1}${\rm)}, then
${\cal A}$ is cyclicly amenable.
\end{list}
\end{theorem}


\vspace{3mm}

\noindent {\sc Hossein Javanshiri}\\
Department of Mathematics,
Yazd University,
P.O. Box: 89195-741, Yazd, Iran.\\
E-mail: h.javanshiri@yazd.ac.ir\\

\noindent {\sc Mehdi Nemati}\\
Department of Mathematical Sciences,
Isfahan Uinversity of Technology,
Isfahan 84156-83111, Iran.\\
E-mail: m.nemati@cc.iut.ac.ir


\footnotesize


\begin{thebibliography}{10}




\bibitem{bad2}
{\sc W. G. Bade and P. C. Curtis},   The Wedderburn decomposition of commutative Banach algebras,
{\it Amer. J. Math.} {\bf 82} (1960), 851--866.


\bibitem{bd}
{W. G. Bade and H. G. Dales}, {The Wedderburn decomposability of
some commutative Banach algebras}, {\it J. Funct. Anal.} {\bf 107}
(1992), 105--121.


\bibitem{bdl}
{W. G. Bade, H. G. Dales and Z. A. Lykova}, {Algebraic and strong
splittings of extensions of Banach algebras}, {\it Mem. Amer.
Math. Soc.} {\bf 137} (1999), 113 pp.



\bibitem{ob}
{O. Berndt}, {On semidirect products of commutative Banach
algebras}, {\it Quaestiones Math.} {\bf 17} (1994), 67--81.



\bibitem{dghg}
{\sc H. G. Dales,  F. Ghahramani and  N. Gr{\o}nb{\ae}k}, Derivations into iterated duals of Banach
algebras, {\it Studia Math.} {\bf 128} (1998), 19--53.


\bibitem{ring0}
{\sc M. D'Anna, C. A. Finocchiaro and M. Fontana}, { Properties of
chains of prime ideals in an amalgamated algebra along an ideal,}
{\it J. Pure Appl. Algebra} {\bf 214} (2010) 1633--1641.


\bibitem{ring1}
{\sc M. D'Anna and M. Fontana}, The amalgamated duplication of a ring along a
multiplicative-canonical ideal, {\it Ark. Mat.} {\bf 45} (2007), 241--252.


\bibitem{ring2}
{\sc M. D'Anna and M. Fontana}, An amalgamated duplication of a ring along an ideal:
The basic properties, {\it J. Algebra Appl.} {\bf 6} (2007), 443--459.


\bibitem{vishki}
{\sc H. R. Ebrahimi Vishki and A. R. Khoddami,} {Character inner
amenability of certain Banach algebras,} {\it Colloq. Math.} {\bf
122} (2011), 225--232.



\bibitem{esfi}
{\sc M. Eshaghi Gordji and M. Filali}, Arens regularity of module actions, {\it Studia math.}
{\bf 181} (2007), 237--254.


\bibitem{fel}
{C. Feldman},  {The Wedderburn principal theorem in Banach
algebras}, {\it Proc. Amer. Math. Soc.} {\bf 2} (1951), 771--777.



\bibitem{ring3}
{\sc R. Fossum}, Commutative extensions by canonical modules are Gorenstein rings,
{\it Proc. Amer. Math. Soc.} {\bf 40} (1973), 395--400.



\bibitem{gro} {\sc N. Gr{\o}nb{\ae}k,} {Weak and cyclic amenability for non-commutative
Banach algebras.} {\it Proc. Edinb. Math. Soc.} {\bf 35} (1992), 315--328.




\bibitem{kani}
 {\sc E. Kaniuth}, A Course in Commutative Banach Algebras,
Graduate Texts in Mathematics, Springer, New York, 2009


\bibitem{kaniuth}
 {\sc E. Kaniuth}, {The Bochner-Schoenberg-Eberlein property and
spectral synthesis for certain Banach algebra products,} {\it
Canad. J. Math.} {\bf 67} (2015), 827--847.


\bibitem{Lau} {\sc A. T.-M. Lau,} {Analysis on a class of Banach algebras
with application to Harmonic analysis on locally compact groups
and semigroups,} {\it Fund. Math.}  {\bf 118} (1983), 161--175.




\bibitem{lauul}
{\sc A. T.-M. Lau and A. \"{U}lger},  Topological centres of certain dual algebras.
{\it Trans. Amer. Math. Soc.} {\bf 348} (1996), 1191--1212.


\bibitem{M} {\sc M. S. Monfared,}
{ On certain products of Banach algebras with applications to
harmonic analysis,} {\it Studia Math.} {\bf 178} (2007), 277--294.


\bibitem{Mon2} {\sc M. S. Monfared,}
{Character amenability of Banach algebras,} {\it Math. Proc.
Cambridge Philos. Soc.} {\bf 144} (2008), 697--706.




\bibitem{nemati1}
{\sc M. Nemati,} { On $\phi$-ergodic property of Banach modules,}
{\it Bull. Belg. Math. Soc. Simon Stevin} {\bf 22} (2015),
655--668.



\bibitem{thomas}
{M. P. Thomas}, {Principal ideals and semi-direct products in
commutative Banach algebras}, {\it J. Funct. Anal.} {\bf 101}
(1991), 312--328.



\bibitem{zhang}
{\sc Y. Zhang}, Weak amenability of module extensions of Banach
algebras, {\it Trans. Amer. Math. Soc.} {\bf 354} (2002),
4131--4151.



\end{thebibliography}
\end{document}